\documentclass[12pt]{article}

\usepackage{amssymb}
\usepackage{graphicx}
\usepackage{epstopdf}
\usepackage{caption}
\newtheorem{theorem}{Theorem}
\newtheorem{lemma}{Lemma}
\newtheorem{remark}{Remark}

\newtheorem{corollary}{Corollary}

\begin{document}
\begin{center}{\bf LIMIT THEOREMS FOR WEIGHTED FUNCTIONALS OF CYCLICAL LONG-RANGE DEPENDENT RANDOM FIELDS\footnote{Short title: Limit Theorems for Cyclical Random Fields}}
\end{center}

\centerline{\bf ANDRIY OLENKO}
\begin{center}{Department of Mathematics and Statistics,
La Trobe University,\\ Victoria, 3086, Australia \\
a.olenko@latrobe.edu.au}
\end{center}

\

 This is an Author's Accepted Manuscript of an article published in the
Stochastic Analysis and Applications, Vol. 31, No. 2. (2013), 199--213.  [copyright Taylor \& Francis], available online at: http://www.tandfonline.com/ [DOI:10.1080/07362994.2013.741410]

\

\begin{center}{\it This contribution is dedicated to the memory of Professor Lakshmikantham.\\ It was organized and communicated by Vo Anh, Member of JSAA.}\end{center}

\

The paper investigates isotropic random fields for which the spectral density is unbounded at some frequencies. Limit theorems for weighted functionals of these random fields are established. It is shown that for a wide class of functionals, which includes the Donsker scheme, the limit is not affected by singularities at non-zero frequencies. For general schemes, in contrast to the Donsker line, we demonstrate that the singularities at non-zero frequencies play a role even for linear functionals.

\

\noindent Keywords: Random fields; limit theorems; long-range dependence; seasonal/cyclical long memory; weighted functionals.

\

\noindent AMS Subject Classification: 60G60, 60F17

\

\section{Introduction}	

Long-range dependence is a well-established empirical fact,
which appears in various fields (finance, signal processing, physics, telecommunications, hydrology,
  etc.), see the monographs \cite{ber,dou,pal} and the numerous references therein.

For a stationary finite-variance random process $\xi(x),\ x\in \mathbb R,$ the most frequently admitted definition of long-range dependence is its non-integrable covariance function $\mathsf B(r)={\rm cov}(\xi(x+r),\xi(x)):$

 \begin{equation}\label{long}\int_{0}^\infty |\mathsf B(r)|\,dr=+\infty\,. \end{equation}
This is often wrongly thought to be due to a singularity of the spectral density at zero frequency.
However,  singularities of the spectral density at non-zero frequencies also imply (\ref{long}).
In particular, models with singularities at non-zero frequencies are of great importance in time series analysis. Many time series show cyclical/seasonal evolutions. It produces peaks in the spectral density whose locations define periods of the cycles, see \cite{anh,rob}.

In spatial statistics models with long-range dependence  have been used to describe a vast number of physical, geological systems and images, see \cite{anh0,fri,gne,has,ma} and references therein. Popular isotropic spatial models with singularities of the spectral density at non-zero frequencies are wave and $J$-Bessel models, see \cite{chi}.

Among the extensive literature on long-range dependence, relatively few publications are devoted to cyclical long-memory processes. Only few cyclical/seasonal long-memory models (GARMA, ARFIMA, ARFISMA, ARUMA) have appeared in the literature.

The case when the summands are some functionals of a long-range dependent Gaussian process is of great importance in the theory of limit theorems for sums of dependent random variables.  It was shown that, compare with Donsker-Prohorov results, long-range dependent summands can produce normalizing coefficients different from $n^{-1/2}$ and non-Gaussian limits.  The non-central limit theorem for functionals of long-range dependent Gaussian processes was investigated by M. Rosenblatt \cite{ros1}, Dobrushin and Major \cite{dob1,dob}, Taqqu \cite{taq1,taq2}, Giraitis and Surgailis \cite{gir1,gir2},  Oppenheim, Haye, and Viano  \cite{ope}. Two volumes by  Doukhan,  Oppenheim, and Taqqu \cite{dou} and by Peccati and Taqqu \cite{pec} give outstanding surveys of the field. For multidimensional results of this type see \cite{leo1,leo2,yad}. 

All mentioned publications were focused on the Donsker line
$$\int_{0}^{rt}H_k(\xi(x))dx,\quad t\in[0,1],$$  where $H_k(\cdot)$ is the $k$-th Hermite polynomial with leading coefficient 1.

In the classical linear case $\int_{0}^{rt}\xi(x)dx,$ $H_1(x)=x,$ the limit is not affected by cyclicity at all. In this case Davydov's theorem \cite{dav} establishes that the limit only depends on the behaviour of the spectral density at
the zero frequency. However, two decades later it was shown that for non-linear functionals with $H_k(\cdot),\ k\ge 2$ the limits can be different  and the cyclical behaviour can play a role, see \cite{dou,haye,ope}.  Some most recent results for stochastic processes can be found  in \cite{iv1} where the asymptotic normality of weighted functionals of Gaussian stationary random processes is obtained. 

The approach taken in the paper continues this line of investigations. The aim of this paper is to show that, even for linear functionals, the limit can be affected by cyclicity. We investigate the limit behaviour of weighted linear functionals of Gaussian random field, when the weights are non-random functions. The Donsker line is a particular case of these schemes.

The paper establishes a limit theorem for random fields, which is also new for the case of stochastic processes and complements one-dimensional results in \cite{iv1}. The limit theorem is obtained under simple conditions on the weight functions.

Several corollaries to the theorem show that
\begin{enumerate}
  \item the limit is not affected by cyclicity for a wide class of functionals, which includes the Donsker scheme;
  \item for general schemes, in contrast to the Donsker line, the cyclical effects play a role even for the linear case $H_1(x)=x.$
\end{enumerate}

The rest of the paper is organized as follows. We begin by describing a class of isotropic random fields with spectral singularities at non-zero frequencies.
 In Section 3, we establish  main results. Finally, further discussion, corollaries and  examples are presented in Section 4.

In what follows we use the symbol $C$ with subscripts to denote constants which are not important for our discussion. Moreover, the same
symbol $C_i$ may be used for different constants appearing in the same proof.

\section{Singularities in spectral densities}
Let $\mathbb{R}^n$ be a Euclidean space of dimension $n \geq 1.$ Furthermore, let $\xi(x),\,x \in \mathbb R^n $ be a real-valued
measurable mean-square continuous homogeneous isotropic Gaussian random field (see \cite{leo2,yad}), with zero mean and the correlation function
 \[\mathsf B_n
(r) = \mathsf B_n (|x|) = \mathsf E \xi (0) \xi (x),\quad   x \in \mathbb
R^n.\]
It is known (see, for example, \cite{leo2,yad}) that there is a bounded nondecreasing function
$\mathbf{\Phi}(u), u \geq 0,$ such that
\[
\mathsf B_n (r) = 2^{\frac{n-2}{2}} \Gamma \Bigl(\frac{n}{2}\Bigr)
\int_0^{\infty} \frac{J_\frac{n-2}{2} (ru)}{(ru)^{\frac{n-2}{2}}}
d \mathbf{\Phi} (u),
\]
where $J_\nu (\cdot) $ is the Bessel function of the first kind of order $\nu \ge
-\frac{1}{2}.$ The function
$\mathbf{\Phi}(\cdot)$ is called the spectral function of the field $\xi(x), x \in \mathbb R^n.$ If there is a function $\varphi(u),\,u\in[0;+\infty)$ such that
\[u^{n-1}\varphi(u)\in L_1([0,\infty)),\quad \mathbf{\Phi}(u)=2
\pi^{n/2}/\,\Gamma(n/2)\int_0^u z^{n-1}\varphi(z)dz,\]
then $\varphi(\cdot)$ is called the isotropic spectral density of the field $\xi(x).$

\begin{figure}[hpb]
\begin{center}
\begin{minipage}{4.1cm}
 \includegraphics[width=4.1cm,trim=2.3cm 1cm 1cm 1cm, clip=true]{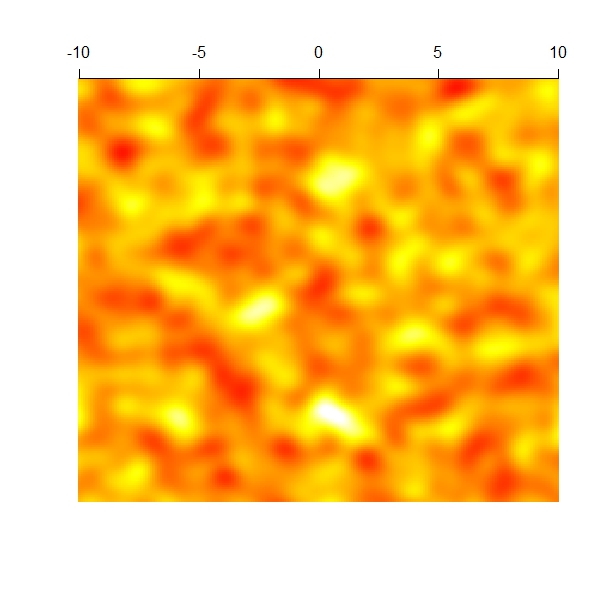}  \end{minipage}
\begin{minipage}{4.1cm}
\includegraphics[width=4.1cm,trim=2.3cm 1cm 1cm 1cm, clip=true]{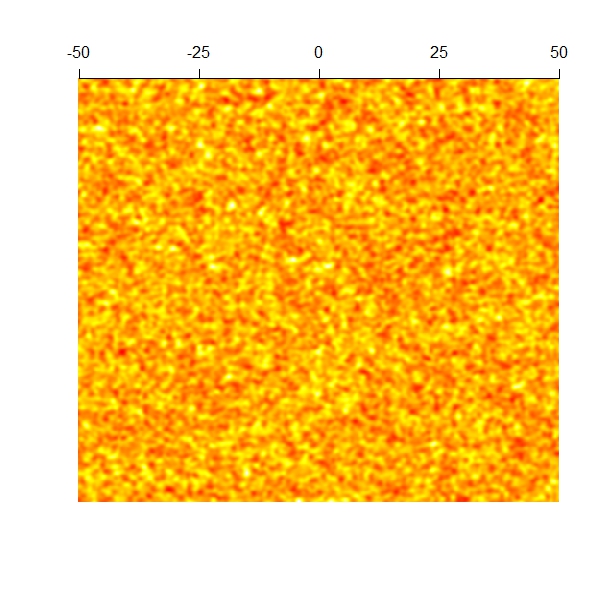} \end{minipage}
\begin{minipage}{4.1cm}
\includegraphics[width=4.1cm,trim=2.3cm 1cm 1cm 1cm, clip=true]{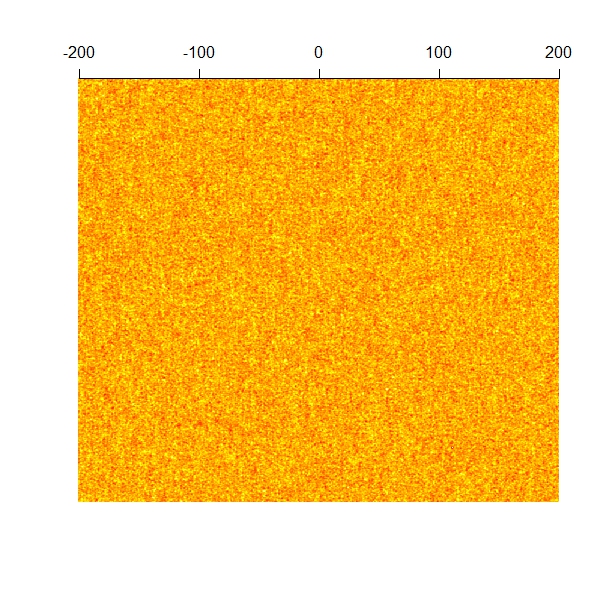} \end{minipage}
\begin{minipage}{4.1cm}
 \includegraphics[width=4.1cm,trim=2.3cm 1cm 1cm 1cm, clip=true]{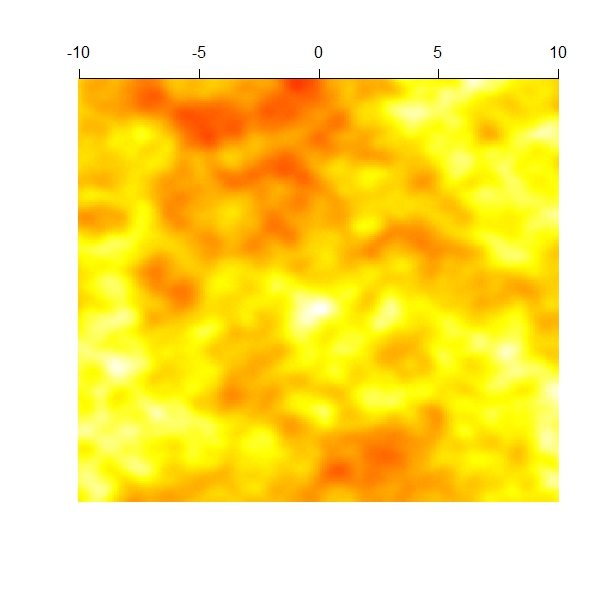}  \end{minipage}
\begin{minipage}{4.1cm}
\includegraphics[width=4.1cm,trim=2.3cm 1cm 1cm 1cm, clip=true]{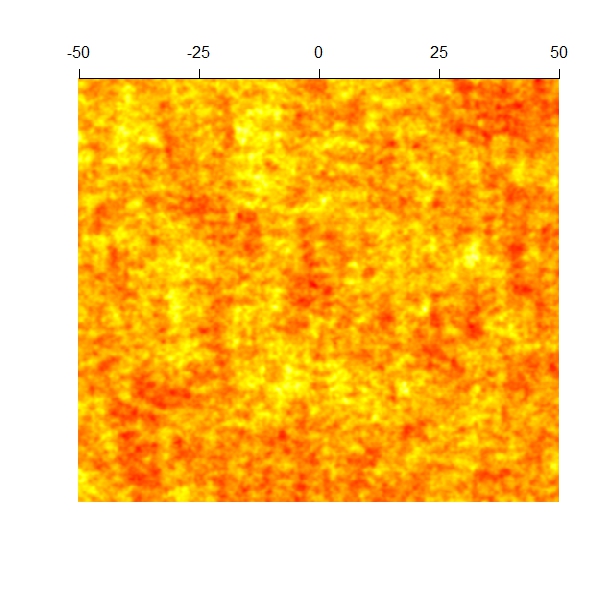} \end{minipage}
\begin{minipage}{4.1cm}
\includegraphics[width=4.1cm,trim=2.3cm 1cm 1cm 1cm, clip=true]{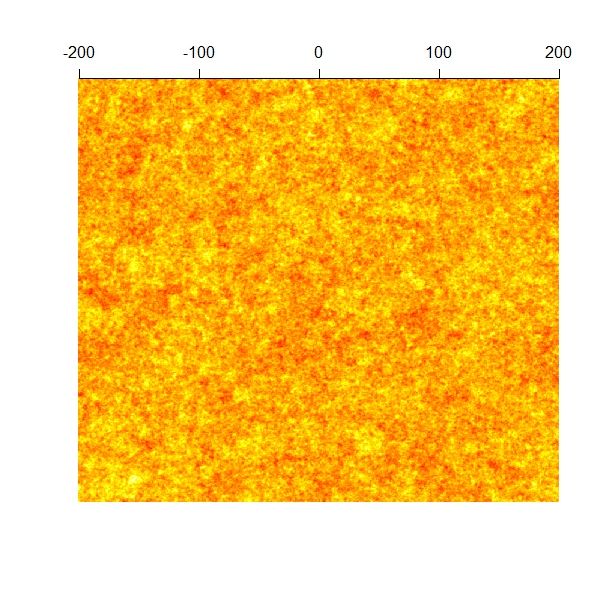} \end{minipage}
\begin{minipage}{4.1cm}
 \includegraphics[width=4.1cm,trim=2.3cm 1cm 1cm 1cm, clip=true]{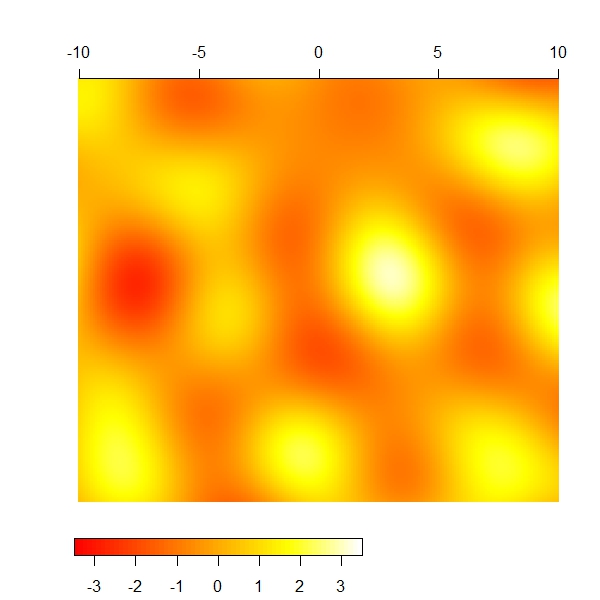}  \end{minipage}
\begin{minipage}{4.1cm}
\includegraphics[width=4.1cm,trim=2.3cm 1cm 1cm 1cm, clip=true]{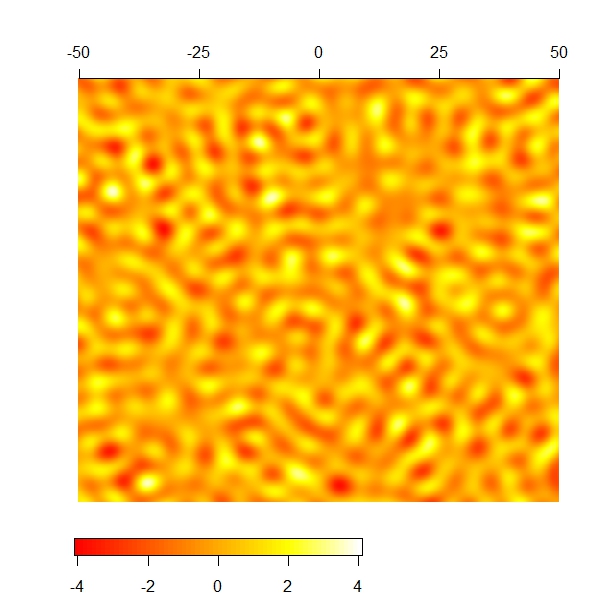} \end{minipage}
\begin{minipage}{4.1cm}
\includegraphics[width=4.1cm,trim=2.3cm 1cm 1cm 1cm, clip=true]{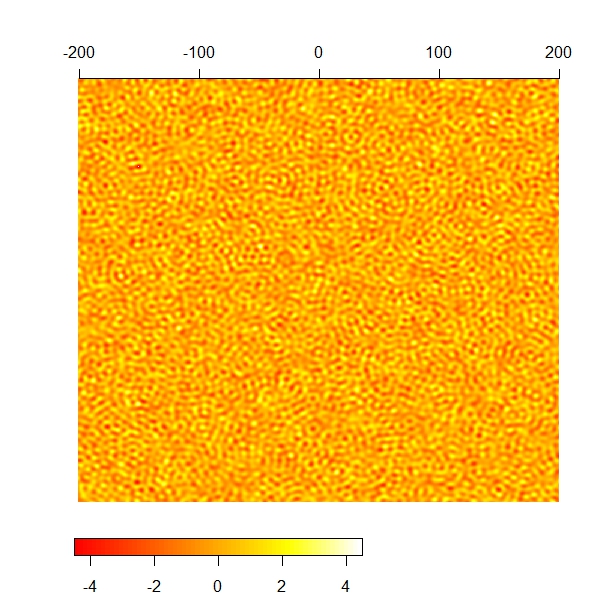} \end{minipage}
\end{center}
  \caption{Two-dimensional realizations of random fields: short-range dependent (first row), long-range dependent (second row), and cyclical long-range dependent (third row)  models}\label{fig1}
\end{figure}

Figure~\ref{fig1} shows two-dimensional realizations of three types of random fields (from top to bottom): short-range dependent normal scale mixture model with nonsingular spectrum, long-range dependent Cauchy model with the spectral singularity at zero, and cyclical long-range dependent wave model with nonzero singularity, see \cite{chi,gne,sch}.  The figure has been generated by the R package \textsc{RandomFields} \cite{sch}.
Images in each row correspond to the same model, but for decreasing distance scales. Each column has the same legend strip for the color scale which is shown at the bottom.  The cyclical behaviour clearly manifests itself in two small scale images in the third row. Such long-range dependent random fields can be used for texture modeling.

Various properties of long-range dependent random fields were investigated in  \cite{leo1,leo2}.  Considered random fields had spectral singularities at zero and the spectral densities of form
  \begin{equation}\label{g1}
     {\varphi}(|\lambda{}|)=\frac{h_0(|\lambda{}|)}{|\lambda{}|^{n-\alpha{}}},
  \end{equation}
where $h_0(\cdot)$ is a bounded function defined on $\mathbb{R}^+:=[0,+\infty),$ $h_0(\cdot)$ is continuous in some neighborhood of
$0,$ and $h_0(0)\ne 0.$

\noindent

We need some analogous of the representation (\ref{g1}) for random fields with singularities of $\mathbf{\Phi}(\lambda{}) $ at the points $a_0,\dots,a_k\in [0,+\infty),$ $a_0=0,\ a_i<a_j,$ for $i<j.$

\

\noindent{\bf Assumption A.}\ \ {\it Let $\xi(x)$ be a Gaussian random field with the spectral density of form
\begin{equation}\label{g5}
    \varphi(|\lambda{}|)=\frac{h_0(|\lambda{}|)}{|\lambda{}|^{n-\alpha{}_0}}+\sum_{i=1}^k\frac{h_i(|\lambda{}|-a_i)}{||\lambda{}|-a_i|^{1-\alpha{}_i}},
\end{equation}
where $\varphi(|\lambda{}|)\in L_1(\mathbb{R}^n),$ $0<\alpha{}_0<n,$ $0<\alpha{}_i<1,$ for $i=1,..,k.$ $h_i(\cdot)\ge 0,\;\;i=0,...,k,$ are bounded functions defined on $[-a_i,+\infty).$ Each $h_i(\cdot)$ is continuous in some neighborhoods of
 $0,$ and $h_i(0)\ne 0.$}

 \begin{remark}
 In this paper we only investigate the case, when the functions $h_i(\cdot)$ are continuous in some neighborhood of $0$ and the power
  $\alpha{}_i$  in the representation {\rm(\ref{g5})} is the same for $|\lambda{}|>a_i$ and $|\lambda{}|<a_i.$
  However, all obtained results can be easily generalized to the case of different left and right limits of the function
   $h_i(\cdot)$ at zero, and  different powers $\alpha{}_i$ for the regions $|\lambda{}|>a_i$ and $|\lambda{}|<a_i.$
 \end{remark}

\begin{remark}\label{rem2}
  We have distinction of kind between spectral singularities of random processes $(n=1)$ and random fields $(n\ge 2).$
   For random processes the spectral density has singularities at $k$ points. For random fields we have a point singularity at $|\lambda|=a_0=0.$ When $|\lambda|=a_i\ne 0$ the spectral density has singularities at all points of the $n$-dimensional sphere $S_n(a_i)=\{\lambda\in\mathbb{R}^n :|\lambda|=a_i\}.$
 \end{remark}

 \begin{remark}
   In contrast to the case of long-range dependent random fields, all powers in the representation
 {\rm(\ref{g5}),} which correspond to nonzero frequency components,  are $1-\alpha{}_i.$
   It follows from Remark~{\rm\ref{rem2}} and the integrability of the spectral density in {\rm(\ref{g5})}.
    \end{remark}

\section{Limit theorem for general schemes}
In this section we investigate the limit behaviour of the weighted  functionals 
$$\int_{\mathbb{R}^n}f_{rt}(x)\xi(x)dx,\quad r\to \infty.$$

We consider radial essentially non-zero weight functions $f_r(\cdot)$ and $f_{j,r}(\cdot),$ i.e.
\begin{itemize}
  \item there are such functions $\tilde{f}_r(\cdot), \tilde{f}_{j,r}(\cdot)$ that $f_r(x)=\tilde{f}_r(|x|),$ $f_{j,r}(x)=\tilde{f}_{j,r}(|x|);$
  \item for each weight function there exists a set of positive measure, where the function is not zero.
\end{itemize}

Let the functions $f_{j,r}:\mathbb{R}^n\to \mathbb{R},$ $j=0,...,k,$ satisfy the conditions:

   \begin{enumerate}
     \item $f_{j,r}(\cdot)\in L_m(\mathbb{R}^n),$ $m=1,2\,;$
     \item the Fourier transform of $r^{-n}f_{j,r}(\cdot)$ is a function $g_{j}(r(|\lambda{}|-a_j)),$ such that
\item $g_j(s), s\in \mathbb{R}$ is an even function, and there exists  $s_0\ge 0$ such that
\begin{equation}
    \label{g}\displaystyle g_j^2(s)\le \frac{C}{s^n},\quad \mbox{for all} \ s\ge s_0.
\end{equation}

   \end{enumerate}
 \begin{lemma}\label{a1} Let the function $f_r:\mathbb{R}^n\to \mathbb{R}$ is given by
\begin{equation}
    \label{f} f_r(x)=\sum_{j=0}^k C_j f_{j,r}(x).
\end{equation}
Then the representation {\rm(\ref{f})} is unique.
 \end{lemma}

{\it Proof.} If there are two different representations in the form (\ref{f}), then, for some $i\in\{0,...,k\},$
\begin{equation}
    \label{fi} f_{i,r}(x)=\sum_{\stackrel{j=0}{j\not= i}}^k \tilde{C}_j f_{j,r}(x),\quad \mbox{for all}\ x\in \mathbb{R}^n.
\end{equation}
Applying the Fourier transform to both sides of (\ref{fi}) we obtain
\begin{equation}
    \label{gi} g_{i}(r(|\lambda|-a_i))=\sum_{\stackrel{j=0}{j\not= i}}^k \tilde{C}_j g_{j}(r(|\lambda|-a_j)),\quad \mbox{for all}\ \lambda\in \mathbb{R}^n.
\end{equation}
For given $s\ge 0$ and $r>0,$ let us choose $|\lambda|=a_i+\frac{s}{r}.$ With this choice of $|\lambda|$ we may rewrite (\ref{gi}) as
$$g_{i}(s)=\sum_{\stackrel{j=0}{j\not= i}}^k \tilde{C}_j g_{j}(s+(a_i-a_j)r).$$
The assumption (\ref{g}) implies that $g_i(s)\equiv 0,$ which is impossible as $f_{i,r}(\cdot)$ is an essentially non-zero function.

Due to Lemma~\ref{a1} we investigate the limit behavior of the functionals in which the weight function $f_{i,r}(x)$ matches the singularity (otherwise the asymptotic distributions are degenerate). Limit theorems for general functionals $\int_{\mathbb{R}^n}f_{rt^{1/n}}\left(x\right)\xi(x)dx$ can be obtained by "factoring" $f_{r}\left(x\right)$ in the form~(\ref{f}) and by applying Lemma~\ref{a1} and Theorem~\ref{a2}.

 \begin{theorem}\label{a2}
Let the isotropic spectral density $\varphi(\cdot)$ satisfy assumption {\rm A}. Then the finite dimensional distributions of the process
   $$X_{r,j}(t)= \left\{
               \begin{array}{ll}
{\displaystyle\frac{1}{\sqrt{A_j\,h_j(0)}r^{n-\alpha{}_j/2}}\int_{\mathbb{R}^n}f_{j,rt^{1/n}}\left(x\right)\xi(x)dx,} & \ t\in(0,1];\vspace{1mm} \\
                 0, &\ t=0
               \end{array}
             \right. $$
 converge weakly to the finite dimensional distributions of the process
   $$ X_j(t)=\left\{
               \begin{array}{ll}
{\displaystyle{t}\int_{\mathbb{R}^n}\frac{g_j(|u|t^{\frac 1 n})}{|u|^{\frac {n-\alpha{}_j} 2}}dZ(u)},& \ t\in(0,1];\vspace{1mm} \\
                 0, &\ t=0\,,
               \end{array}
             \right. $$
when $r\to\infty.$ $Z(\cdot)$ is the Wiener measure on $(\mathbb{R}^n, \mathfrak{B}^n),$ $$A_j=\left\{
               \begin{array}{ll}
               1, &\ j=0,\\
                 {2a_j^{n-1}}, &\ j\not=0.
               \end{array}
             \right. $$
 \end{theorem}

{\it Proof.} 
We first use the spectral representation of random fields
\begin{equation}
    \label{g8} \xi(x)=\int_{\mathbb{R}^n}e^{i<\lambda{},x>}\sqrt{\varphi(|\lambda{}|)}dW(\lambda{}),
\end{equation}
where $W(\cdot) $ is the Wiener measure on  $(\mathbb{R}^n, \mathfrak{B}^n),$ see \cite{leo2}.

 Let us define
\begin{equation} \label{g9}
I_j(r):=\left\{
               \begin{array}{ll}
{\displaystyle{r^{-n}}\int_{\mathbb{R}^n}f_{j,r}\left(x\right)\xi(x)dx}, & \ r>0;\vspace{1mm} \\
                 0, &\ r=0\,.
               \end{array}
             \right.
\end{equation}

To investigate the limit behaviour of $I_j(rt^\frac 1n)$  we generalize the approach of \cite{olkly} to the case of arbitrary functions $g_j$ and spectral singularities at nonzero frequencies.

Heuristically, by substituting the spectral representation (\ref{g8}) in (\ref{g9}) and by changing the order of integration, we obtain

 $$     r^{-n}\int_{\mathbb{R}^n}f_{j,r}\left(x\right)\int_{\mathbb{R}^n}e^{i<\lambda{},x>}\sqrt{\varphi(|\lambda{}|)}dW(\lambda{})dx$$
\begin{equation}  \label{g10}     \stackrel{d}{=}r^{-n}\int_{\mathbb{R}^n}\sqrt{\varphi(|\lambda{}|)}\int_{\mathbb{R}^n}e^{i<\lambda{},x>}f_{j,r}\left(x\right)dx dW(\lambda{}).
\end{equation}

To prove  (\ref{g10}) we note, that
  $$ \int_{\mathbb{R}^n}e^{i<\lambda{},x>}\sqrt{\varphi(|\lambda{}|)}dW(\lambda{}) $$
is the Fourier transform of the stochastic measure $\mu(d\lambda{})=\sqrt{\varphi(|\lambda{}|)}dW(\lambda{}),$
and $$\int_{\mathbb{R}^n}e^{i<\lambda{},x>}f_{j,r}\left(x\right)dx $$ is the Fourier transform of the non-random function $f_{j,r}\left(x\right).$

Note that the formula (\ref{g10}) can be rewritten as
$$\int_{\mathbb{R}^n}f_{j,r}\left(x\right)d\widehat{\mu}(x)\stackrel{d}{=}\int_{\mathbb{R}^n}\widehat{f}_{j,r}(\lambda{})d\mu(\lambda{}). $$
This is a stochastic analogous of Plancherel's  formula.  Since $f_{j,r}(\cdot)\in L_2(\mathbb{R}^n)$ and
$\sqrt{\varphi(|\cdot|)}\in L_2(\mathbb{R}^n),$ (\ref{g10}) is a correct transformation, see \cite{hou}.

The function  $g_j
(r(|\lambda{}|-a_j))$ is the Fourier transform of $r^{-n}f_{j,r}(x).$ Therefore we get
\begin{equation} \label{g12j}
      r^{-n} \int_{\mathbb{R}^n}\int_{\mathbb{R}^n}e^{i<\lambda{},x>}f_{j,r}\left({x}\right)dx\sqrt{\varphi(|\lambda{}|)}dW(\lambda{})
       \stackrel{d}{=}\int_{\mathbb{R}^n}g_j
(r(|\lambda{}|-a_j))\sqrt{\varphi(|\lambda{}|)}dW(\lambda{}).
\end{equation}

The integral on the right-hand  side of (\ref{g12j}) can be partitioned as follows
$$ \int_{\mathbb{R}^n} g_j(r(|\lambda{}|-a_j))\sqrt{\varphi(|\lambda{}|)}W(d\lambda{})
        \stackrel{d}{=} \int_{|\lambda{}|>a_j} g_j(r(|\lambda{}|-a_j))\sqrt{\varphi(|\lambda{}|)}dW(\lambda{})$$
$$+ \int_{|\lambda{}|<a_j} g_j(r(|\lambda{}|-a_j))\sqrt{\varphi(|\lambda{}|)}dW(\lambda{})=:K_1(r)+K_2(r).$$

First we deal with $K_1(r).$
By making the change of variables
\begin{equation} \label{g13}\left\{%
      \begin{array}{ll}
         u=\lambda{}\left(1-\frac{a_j}{|\lambda{}|}\right), & |u|=|\lambda{}|-a_j \\
         \lambda{}=u\left( 1+\frac{a_j}{|u|}\right), & |\lambda{}|=|u|+a_j \\
      \end{array}%
           \right.
\end{equation}
we bijectively map the set
$\{ \lambda{}\in\mathbb{R}^n: |\lambda{}|>a_j\}$
into the set
$\mathbb{R}^n \backslash\{(0,...,0)\}.$

\begin{lemma}{\rm\cite{olkly}}
The Jacobian of the transformation {\rm(\ref{g13})} is
    $$ \det(J_n(u))= \left( 1+\frac{a_j}{|u|}\right)^{n-1},\;\; |u|\neq 0.$$
\end{lemma}

\begin{corollary}\label{nas1} By {\rm(\ref{g5})} and $\varphi(\lambda)\in L_1(\mathbb{R}^n)$ it follows that
$$\int_0^\infty\frac{h_0\left(\rho+a_j\right)}
{\left(\rho+a_j\right)^{n-\alpha{}_0}}\left(1+\frac{a_j}{\rho}\right)^{n-1}\rho^{n-1}d\rho<+\infty\,;$$

$$\int_0^\infty\frac{h_i\left(\rho-a_i+a_j\right)}{\left|\rho-a_i+a_j\right|^{1-\alpha{}_i}}\left(1+\frac{a_j}{\rho}\right)^{n-1}\rho^{n-1}d\rho<+\infty\,.$$
\end{corollary}

Let us use the change of variable formula for stochastic integrals, see \cite[Proposition 4.2]{dob1} and \cite[Theorem 4.4]{maj}. Note, that the integrand involves only function in our case. We get
\begin{eqnarray*}
       K_1(r)& \stackrel{d}{=}& \int_{\mathbb{R}^n} \sqrt{\left(\frac{h_0(|u|+a_j)}{\left(|u|+a_j\right)^{n-\alpha{}_0}}+\sum_{i=1}^k\frac{h_i(|u|-a_i+a_j)}{||u|-a_i+a_j|^{1-\alpha{}_i}}\right)}\\
       && \times g_j (r|u|)\left(1+\frac{a_j}{|u|}\right)^{\frac{n-1}{2}}d\widehat{W}(u),
    \end{eqnarray*}
where $\widehat{W}(\cdot)$ is the Wiener measure on $(\mathbb{R}^n, \mathfrak{B}^n).$

By making the change of variables  $ru_i=\widetilde{u}_i,\;\;i=\overline{1,n}$ in the last integral and by the self-similarity of Gaussian white noise,
$d\widehat{W}(cx)\stackrel{d}{=}c^\frac{n}{2} d\widetilde{W}(x)\,,$
we obtain
$$K_1(rt^\frac 1n ) \stackrel{d}{=} \int_{\mathbb{R}^n} r^{-\frac{\alpha{}_j}{2}} \frac{g_j(|\widetilde{u}|t^\frac 1n)}{ |\widetilde{u}|^\frac {n-\alpha_j}{2}}\left(I_{\{j\not= 0\}}\,h_j\left(\frac{|\widetilde{u}|}{r}\right)+\left(\frac{|\widetilde{u}|}{r}\right)^{1-\alpha{}_j}\frac{h_0\left(\frac{|\widetilde{u}|}{r}+a_j\right)}{\left(\frac{|\widetilde{u}|}{r}+a_j\right)^{n-\alpha{}_0}}\right.
$$
\begin{equation} \label{g15b}\left.+\mathop{\sum_{i=1}^k}_{i\not=j} \left(\frac{|\widetilde{u}|}{r}\right)^{1-\alpha_j} \frac{h_i\left(\frac{|\widetilde{u}|}{r}-a_i+a_j\right)}{\left|\frac{|\widetilde{u}|}{r}-a_i+a_j\right|^{1-\alpha{}_i}}\right)^{1/2}
\left(a_j+\frac{|\widetilde{u}|}{r}\right)^{\frac{n-1}{2}}d\widetilde{W}(\widetilde{u}),
\end{equation}
where
$$I_{\{j\not= 0\}}=
\left\{%
      \begin{array}{ll}
         1, \; & \mbox{if}\;\; j\not= 0, \\
         0,\; & \mbox{if} \;\;j = 0. \\
      \end{array}%
           \right.
  $$

Denote
   \begin{equation} \label{g16b}
       {Y}_{j,1}(t):=\left\{
               \begin{array}{ll}
{\displaystyle\int_{\mathbb{R}^n}\frac{g_j(|\widetilde{u}|t^\frac 1n)}{|\widetilde{u}|^{\frac {n-\alpha_j}{2} }}d\widetilde{W}(\widetilde{u}),}
& \ t\in(0,1];\vspace{1mm} \\
                 0, &\ t=0\,.
               \end{array}
             \right.
   \end{equation}
By (\ref{g15b}) and (\ref{g16b}),
\begin{equation} \label{g17a}
      {R}_r(t):=\mathsf{E}\left( t\sqrt{\frac{r^{\alpha{}_j}}{V_j\,h_j(0)}}K_1(rt^\frac 1n)-{t}{Y}_{j,1}(t)   \right)^2
         = t^2\int_{\mathbb{R}^n}\frac{g^2_j(|\widetilde{u}|t^\frac 1n)}{|\widetilde{u}|^{n- \alpha{}_j}} {Q}_r(|\widetilde{u}|)d\widetilde{u},
\end{equation}
where
$$V_j=\left\{
               \begin{array}{ll}
               1, &\ j=0,\\
                 {a_j^{n-1}}, &\ j\not=0,
               \end{array}
             \right. $$
$${Q}_r(|\widetilde{u}|):=\left({V_j^{-1/2}}\left(a_j+\frac{|\widetilde{u}|}{r}\right)^{\frac{n-1}{2}}\left(I_{\{j\not= 0\}}\,\frac{h_j\left(\frac{|\widetilde{u}|}{r}\right)}{h_j(0)}+\frac{\left(\frac{|\widetilde{u}|}{r}\right)^{1-\alpha_j}h_0\left(\frac{|\widetilde{u}|}{r}+a_j\right)}
{h_j(0)\left(\frac{|\widetilde{u}|}{r}+a_j\right)^{n-\alpha{}_0}}\right.\right.$$
$$\left.\left.+\mathop{\sum_{i=1}^k}_{i\not=j} \left(\frac{|\widetilde{u}|}{r}\right)^{1-\alpha_j} \frac{h_i\left(\frac{|\widetilde{u}|}{r}-a_i+a_j\right)}{h_j(0)\left|\frac{|\widetilde{u}|}{r}-a_i+a_j\right|^{1-\alpha{}_i}}\right)^{1/2}
-1\right)^2.$$

Let us choose $\psi(r)>0$  so that $\psi(r)\to\infty,$ $\frac{\psi(r)}{r}\to 0,$ and $\frac{r^{\alpha_j}}{\psi^n(r)}\to 0,$ when $r\to\infty.$
If $\mathbb{R}^n$ is split into the regions
$$B_1:=\{\widetilde{u}\in\mathbb{R}^n: |\widetilde{u}|\le \psi(r)\} \quad \mbox{and}\quad
B_2:=\{\widetilde{u}\in\mathbb{R}^n: |\widetilde{u}|> \psi(r)\}, $$
then the integral on the right-hand side of (\ref{g17a}) becomes the sum $R_r(t)=R_{r,1}(t)+R_{r,2}(t).$

By assumption A, given any $\varepsilon>0,$ there is  $r_0>0$ such that $Q_r(|\widetilde{u}|)<\varepsilon,$ when $r>r_0,\;\widetilde{u}\in B_1.$

 $L_1$ and $L_2$-integrability of the function $f_{j,r}(\cdot)$ yields  $g_j^2(|\lambda{}|)\in L_2(\mathbb{R}^n)$ and $g_j^2(|\lambda{}|)$ is bounded. Hence
 the integral
$$ t^2\int_{\mathbb{R}^n}\frac{g_j^2
(|\widetilde{u}|t^{1/n})}{|\widetilde{u}|^{n- \alpha{}_j  }} d\widetilde{u}=t^{2-\frac{\alpha{}_j}{n}}\int_{\mathbb{R}^n}\frac{g_j^2
(|\widetilde{z}|)}{|\widetilde{z}|^{n- \alpha{}_j}} d\widetilde{z},
$$
is uniformly bounded on $t\in[0,1].$ Thus
 $R_{r,1}(t)$ can be made arbitrarily small  by decreasing the value of $\varepsilon.$

The decay rate of $g_j^2(\cdot)$ implies that there is $r^0>0 $ such that, for any $r>r^0,$
$${R}_{r,2}(t)\le t^2\int_{B_2}\frac{g_j^2
(|\widetilde{u}|t^{1/n})}{|\widetilde{u}|^{n- \alpha{}_j }} \left(\left(I_{\{j\not= 0\}}\,\frac{h_j\left(\frac{|\widetilde{u}|}{r}\right)}{h_j(0)}+\left(\frac{|\widetilde{u}|}{r}\right)^{1-\alpha_j}\frac{h_0\left(\frac{|\widetilde{u}|}{r}+a_j\right)}
{h_j(0)\left(\frac{|\widetilde{u}|}{r}+a_j\right)^{n-\alpha{}_0}}\right.\right.$$
$$\left.\left.+\mathop{\sum_{i=1}^k}_{i\not=j} \left(\frac{|\widetilde{u}|}{r}\right)^{1-\alpha_j} \frac{h_i\left(\frac{|\widetilde{u}|}{r}-a_i+a_j\right)}{h_j(0)\left|\frac{|\widetilde{u}|}{r}-a_i+a_j\right|^{1-\alpha{}_i}}\right){V_j^{-1}}\left(a_j+\frac{|\widetilde{u}|}{r}\right)^{n-1}+1\right)d\widetilde{u}$$
$$\le \int_{B_2}\frac{C_1dz}{|{z}|^{2n- \alpha{}_j }}+ \frac{C_2}{\psi^n(r)r^{-\alpha_j}}\left(\int_{\frac{\psi(r)}{r}}^\infty I_{\{j\not= 0\}}\,\frac{h_j\left(\rho\right)}{\rho^{1-\alpha_j}}
\left(1+\frac{a_j}{\rho}\right)^{n-1}\rho^{n-1}
d\rho\right.$$
$$\left.+\int_{\frac{\psi(r)}{r}}^\infty\left(\frac{h_0\left(\rho+a_j\right)}
{\left(\rho+a_j\right)^{n-\alpha{}_0}}+\mathop{\sum_{i=1}^k}_{i\not=j}  \frac{h_i\left(\rho-a_i+a_j\right)}{\left|\rho-a_i+a_j\right|^{1-\alpha{}_i}}\right)\left(1+\frac{a_j}{\rho}\right)^{n-1}\rho^{n-1}d\rho\right)\,. $$

Now, by corollary~\ref{nas1} and the choice of $\psi(r)$ the upper bound vanishes when $r\to+\infty.$
Thus, $\lim_{r\to\infty}{R}_r(t)=0$ and the finite dimensional distributions of the process $t\sqrt{\frac{r^{\alpha{}_j}}{h_j(0)V_j}}K_1(rt^\frac 1n)$ converge to the finite dimensional distributions of the process $t\,{Y}_{j,1}(t).$

Since the integral $K_2(\cdot)=0$ if $j=0,$ we  investigate  it only for $j\not= 0.$
By making the change of variables
\begin{equation} \label{g18-} \left\{%
      \begin{array}{ll}
         u=\lambda{}\left(\frac{a_j}{|\lambda{}|}-1\right), & |u|=a_j-|\lambda{}|, \\
         \lambda{}=u\left( \frac{a_j}{|u|}-1\right), & |\lambda{}|=a_j-|u| \\
      \end{array}%
           \right.
\end{equation}
we bijectively map the set $\{\lambda{}\in\mathbb{R}^n: 0<|\lambda{}|<a_j\}$
into itself.

\begin{lemma} {\rm\cite{olkly}}
The Jacobian of the transformation {\rm(\ref{g18-})} is
    $$ \det(\widetilde{J}_n(u))= -\left( \frac{a_j}{|u|}-1\right)^{n-1},\;\; |u|\neq 0.$$
\end{lemma}

\begin{corollary}
$$\int_0^{a_j}\frac{h_0\left(a_j-\rho\right)}
{\left(a_j-\rho\right)^{n-\alpha{}_0}}\left(\frac{a_j}{\rho}-1\right)^{n-1}\rho^{n-1}d\rho<+\infty\,;$$

$$\int_0^{a_j}\frac{h_i\left(a_j-\rho-a_i\right)}{\left|a_j-\rho-a_i\right|^{1-\alpha{}_i}}\left(\frac{a_j}{\rho}-1\right)^{n-1}\rho^{n-1}d\rho<+\infty\,.$$
\end{corollary}

Just as in the case of $K_1(\cdot) $ one can show that
 $$K_2(rt^\frac 1n ) \stackrel{d}{=} \int_{|\widetilde{u}|<ra_j} r^{-\frac{\alpha_j}{2}} a_j^\frac{n-1}{2}\frac{g_j(|\widetilde{u}|t^\frac 1n)}{ |\widetilde{u}|^\frac {n-\alpha_j}{2}} \left(h_j\left(-\frac{|\widetilde{u}|}{r}\right)+\left(\frac{|\widetilde{u}|}{r}\right)^{1-\alpha_j}\frac{h_0\left(a_j-\frac{|\widetilde{u}|}{r}\right)}
 {\left(a_j-\frac{|\widetilde{u}|}{r}\right)^{n-\alpha{}_0}}\right.
 $$
\begin{equation} \label{g15b-}\left.+\mathop{\sum_{i=1}^k}_{i\not=j} \left(\frac{|\widetilde{u}|}{r}\right)^{1-\alpha_j} \frac{h_i\left(a_j-\frac{|\widetilde{u}|}{r}-a_i\right)}{\left|a_j-\frac{|\widetilde{u}|}{r}-a_i\right|^{1-\alpha{}_i}}\right)^{1/2}
\left(1-\frac{|\widetilde{u}|}{ra_j}\right)^{\frac{n-1}{2}}d\overline{W}(\widetilde{u})\,,
\end{equation}
where $\overline{W}(\cdot)$ is the Wiener measure on $(\mathbb{R}^n, \mathfrak{B}^n).$

Denote
   \begin{equation} \label{g16b-}
       Y_{j,2}(t):=\left\{
               \begin{array}{ll}
{\displaystyle \int_{\mathbb{R}^n}\frac{g_j(|\widetilde{u}|t^\frac 1n)}{|\widetilde{u}|^{\frac {n-\alpha{}_j}{2} }}d\overline{W}(\widetilde{u}),}
& \ t\in(0,1];\vspace{1mm} \\
                 0, &\ t=0\,.
               \end{array}
             \right.
   \end{equation}
Combining (\ref{g15b-}) and (\ref{g16b-}), we obtain
\begin{equation} \label{g17b-}
      {S}_r(t):=\mathsf{E}\left( t\sqrt{\frac{r^{\alpha{}_j}}{h_j(0)a_j^{n-1}}}K_2(rt^\frac 1n)-{t}{Y}_{j,2}(t)   \right)^2
         = t^2\int_{\mathbb{R}^n}\frac{g^2_j(|\widetilde{u}|t^\frac 1n)}{|\widetilde{u}|^{n- \alpha{}_j}} \overline{Q}_r(|\widetilde{u}|)d\widetilde{u},
\end{equation}
where
$$\overline{Q}_r(|\widetilde{u}|):=\left(\left(1-\frac{|\widetilde{u}|}{ra_j}\right)_+^{\frac{n-1}{2}}\left(\frac{h_j\left(-\frac{|\widetilde{u}|}{r}\right)}{h_j(0)}+\left(\frac{|\widetilde{u}|}{r}\right)^{1-\alpha_j}\frac{h_0\left(a_j-\frac{|\widetilde{u}|}{r}\right)}
{h_j(0)\left(a_j-\frac{|\widetilde{u}|}{r}\right)^{n-\alpha{}_0}}\right.\right.$$
$$\left.\left.+\mathop{\sum_{i=1}^k}_{i\not=j} \left(\frac{|\widetilde{u}|}{r}\right)^{1-\alpha_j} \frac{h_i\left(a_j-\frac{|\widetilde{u}|}{r}-a_i\right)}{h_j(0)\left|a_j-\frac{|\widetilde{u}|}{r}-a_i\right|^{1-\alpha{}_i}}\right)^{1/2}-1\right)^2\,,$$
$$(x)_+=\left\{%
      \begin{array}{ll}
         x, \; & \mbox{if}\;\; x> 0, \\
         0,\; & \mbox{if} \;\;x\le 0. \\
      \end{array}%
           \right.
  $$

The integral on the right-hand side of (\ref{g17b-}) can be split into two parts
${S}_r(t)={S}_{r,1}(t)+{S}_{r,2}(t),$ where the
integration sets are respectively $B_1$ and $B_2.$    Similarly to the case of $R_{r,1}(t)$  the integral
$S_{r,1}(t)$  can be made arbitrarily small  by decreasing the value of  $\varepsilon.$

By (\ref{g}) there is $r^0>0,$ such that, for any $r>r^0,$
$${S}_{r,2}(t)\le t^2\int_{B_2}\frac{g_j^2
(|\widetilde{u}|t^{1/n})}{|\widetilde{u}|^{n- \alpha{}_j }}
\left(\left(\frac{h_j\left(-\frac{|\widetilde{u}|}{r}\right)}{h_j(0)}+\left(\frac{|\widetilde{u}|}{r}\right)^{1-\alpha_j}\frac{h_0\left(a_j-\frac{|\widetilde{u}|}{r}\right)}
{h_j(0)\left(a_j-\frac{|\widetilde{u}|}{r}\right)^{n-\alpha{}_0}}\right.\right.$$
$$\left.\left.+\mathop{\sum_{i=1}^k}_{i\not=j} \left(\frac{|\widetilde{u}|}{r}\right)^{1-\alpha_j} \frac{h_i\left(a_j-\frac{|\widetilde{u}|}{r}-a_i\right)}{h_j(0)\left|a_j-\frac{|\widetilde{u}|}{r}-a_i\right|^{1-\alpha{}_i}}\right)\left(1-\frac{|\widetilde{u}|}{ra_j}\right)_+^{n-1}+1\right)d\widetilde{u}$$
$$\le \int_{B_2}\frac{C_1dz}{|{z}|^{2n- \alpha{}_j }}+ \frac{C_2}{\psi^n(r)r^{-\alpha_j}}\left(\int_{\frac{\psi(r)}{r}}^{a_j} \frac{h_j\left(-\rho\right)}{\rho^{1-\alpha_j}}
\left(\frac{a_j}{\rho}-1\right)^{n-1}\rho^{n-1}
d\rho\right.$$
$$\left.+\int_{\frac{\psi(r)}{r}}^{a_j}\left(\frac{h_0\left(a_j-\rho\right)}
{\left(a_j-\rho\right)^{n-\alpha{}_0}}+ \mathop{\sum_{i=1}^k}_{i\not=j}  \frac{h_i\left(a_j-\rho-a_i\right)}{\left|a_j-\rho-a_i\right|^{1-\alpha{}_i}}\right)\left(\frac{a_j}{\rho}-1\right)^{n-1}\rho^{n-1}d\rho\right)\,. $$

The upper bound vanishes when $r\to+\infty.$

Therefore, the finite dimensional distributions of
$t\sqrt{\frac{r^{\alpha{}_j}}{h_j(0)a_j^{n-1}}}K_2(rt^\frac 1n)$\
converge to the finite dimensional distributions of the process ${t}\,Y_{j,2}(t),\; t\in[0,1],$ when $r\to+\infty.$

The above results imply the convergence of the finite dimensional distributions of the process
$$ \sqrt{2}\, X_{r,j}(t)\stackrel{d}{=} t\sqrt{\frac{r^{\alpha{}_j}}{h_j(0)V_j}}\left(K_1(rt^\frac 1n)+K_2(rt^\frac 1n)\right) $$
to the finite dimensional distributions of
$${t}  ( Y_{j,1}(t)+I_{\{j\not= 0\}}\,Y_{j,2}(t)), \;\;t\in[0,1].$$

Finally,  since the Wiener measures $\widetilde{W}(\cdot)$ and $\overline{W}(\cdot)$ are independent, it follows that there exists the Wiener measure $Z(\cdot) $ on $(\mathbb{R}^n, \mathfrak{B}^n)$
such that, formally,
$$ \sqrt{2}dZ(\cdot)\stackrel{d}{=} d\widetilde{W}(\cdot)+d\overline{W}(\cdot ).$$

Now, it is easy to see that
\begin{eqnarray*}
       {t}  ( Y_{j,1}(t)+(j)_+\,Y_{j,2}(t)) &= {t}\left( {\displaystyle\int_{\mathbb{R}^n}}\frac{g_j(|\widetilde{u}|t)}{|\widetilde{u}|^{\frac {n-\alpha{}_j} 2}}d\widetilde{W}(\widetilde{u})+ (j)_+\, {\displaystyle\int_{\mathbb{R}^n}}\frac{g_j(|\widetilde{u}|t)}{|\widetilde{u}|^{\frac {n-\alpha{}_j} 2}}d\overline{W}(\widetilde{u})\right) \\ &
        \stackrel{d}{=}
      \sqrt{2}^{I_{\{j\not= 0\}}}\,t {\displaystyle\int_{\mathbb{R}^n}}\frac{g_j(|\widetilde{u}|t)}{|\widetilde{u}|^{\frac {n-\alpha{}_j} 2}}dZ(\widetilde{u})=\sqrt{2}^{I_{\{j\not= 0\}}}\,X_j(t).
\end{eqnarray*}

\section{Corollaries and Discussion}

In the particular case when the weight functions are from $f_{0,r}(\cdot)$ class, then the limit $X(t)$ is not affected by cyclicity at all. It is easy to see, that the weight function $f_{0,r}(x)$ is of the form $f\left(\frac{|x|}{r}\right).$ Two examples of such well-known schemes are given below. In both cases the functions $f_{0,r}(x)$ have finite supports.

 \begin{enumerate}
  \item  The Donsker line is a particular case of our general results with
\begin{equation} \label{f0}
f_{0,r}(x)=\left\{
               \begin{array}{ll}
                 1, & \hbox{if}\ |x|\le r \\
                 0, & \hbox{if}\ |x|>r
               \end{array}
             \right.\quad \mbox{and}\quad g_0
(r(|\lambda{}|))=(2\pi)^{n/2}\frac{J_{\frac{n}{2}}(r(|\lambda{}|))}{(r(|\lambda{}|))^{n/2}}.
\end{equation}
\item Functionals of the form $\int_{\mathbb{R}^n}f(x,r)H_m(\xi(x))dx,$ where $\xi(\cdot)$ has only singularity at zero frequency,
$$f(x,r)= \left\{
                                       \begin{array}{ll}
                                         f\left(\frac{|x|}{r}\right), & \hbox{if}\ x\in \triangle(r); \\
                                         0, & \hbox{if}\ x\not\in \triangle(r);
                                       \end{array}
                                     \right.
$$
and $\{\triangle(r), r>0\}$ is an increasing set of homothetic convex bodies,
were investigated in \cite{leo1}.
For the considered case $m=1$ our class of weight functions is wider and the conditions in our limit theorem are much simpler than the corresponding ones in \cite{leo1}.
\end{enumerate}

A particular example of random fields with only singularity at the frequency $a_1\not= 0$ and the Bessel weight function was studied in \cite{ole1} and \cite{olkly}.

 One can easily obtain new examples of the weight functions $f_{j,r}(x)$ and $g_{j}(r(|\lambda|-a_j)),$  by using the representation of the Fourier transform of radial functions \cite{boch} and tables of the Hankel transforms (for example, \cite{erd}).

The obtained results can be straightforward translated to discrete settings, in particular to one-dimensional time series, with sums instead of  integrals in the definition of $X_{r,j}(t)$ in Theorem~\ref{a2}.

In the particular case of stochastic processes ($n=1$) the obtained results complement \cite{iv1}. The article \cite{iv1} investigated the asymptotic normality of the corresponding weighted functionals of non-linear
transformations of Gaussian stationary random processes. The approach was based on the Central Limit Theorem by Peccati
and Tudor \cite{pec0}. It is worth to mention that some partial results for  $k=2$ and discrete time processes were obtained in \cite{ope}, where it was shown that the limits are non-Gaussian (the Rosenblatt process or sums of two independent Rosenblatt processes). An interesting non-trivial problem is the generalization of approaches in the paper and \cite{iv1} to describe all possible types of limit bahaviour for  dimensions $n>1$ and Hermite polynomial degrees~$k>1.$

It would be interesting to extend these approaches to study the limit behaviour of systems modelled by stochastic differential equations, see \cite{anh0,lad}.

Some potential statistical applications of obtained results are in weighted regression analysis and asymptotic inference of fields with singular spectra, see \cite{anh,iv1,leo1,lav}  and references therein. Two possible scenarios are: to choose weight functions with finite supports which correspond to increasing observation regions and to assign values of weight functions by the relative density of observations at certain subregions.

\begin{figure}[h]
\begin{minipage}{7.1cm}
 \includegraphics[width= 6.1cm,trim=0.6cm 1.2cm 1cm 1.4cm, clip=true]{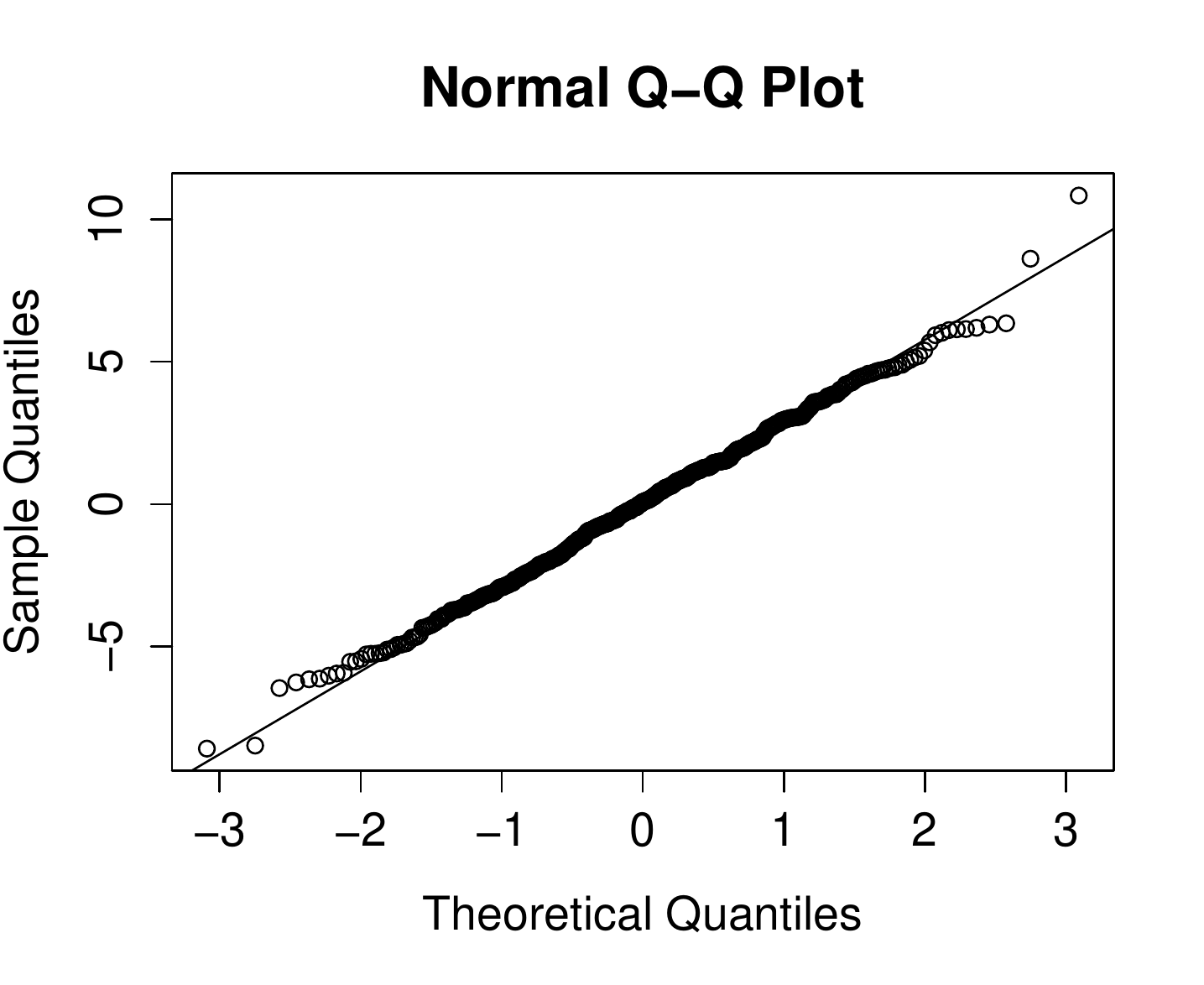}  \end{minipage}\hspace{1cm}
\begin{minipage}{ 7.1cm}
 \includegraphics[width= 6.1cm,trim=1cm 1.2cm 1cm 1.4cm, clip=true]{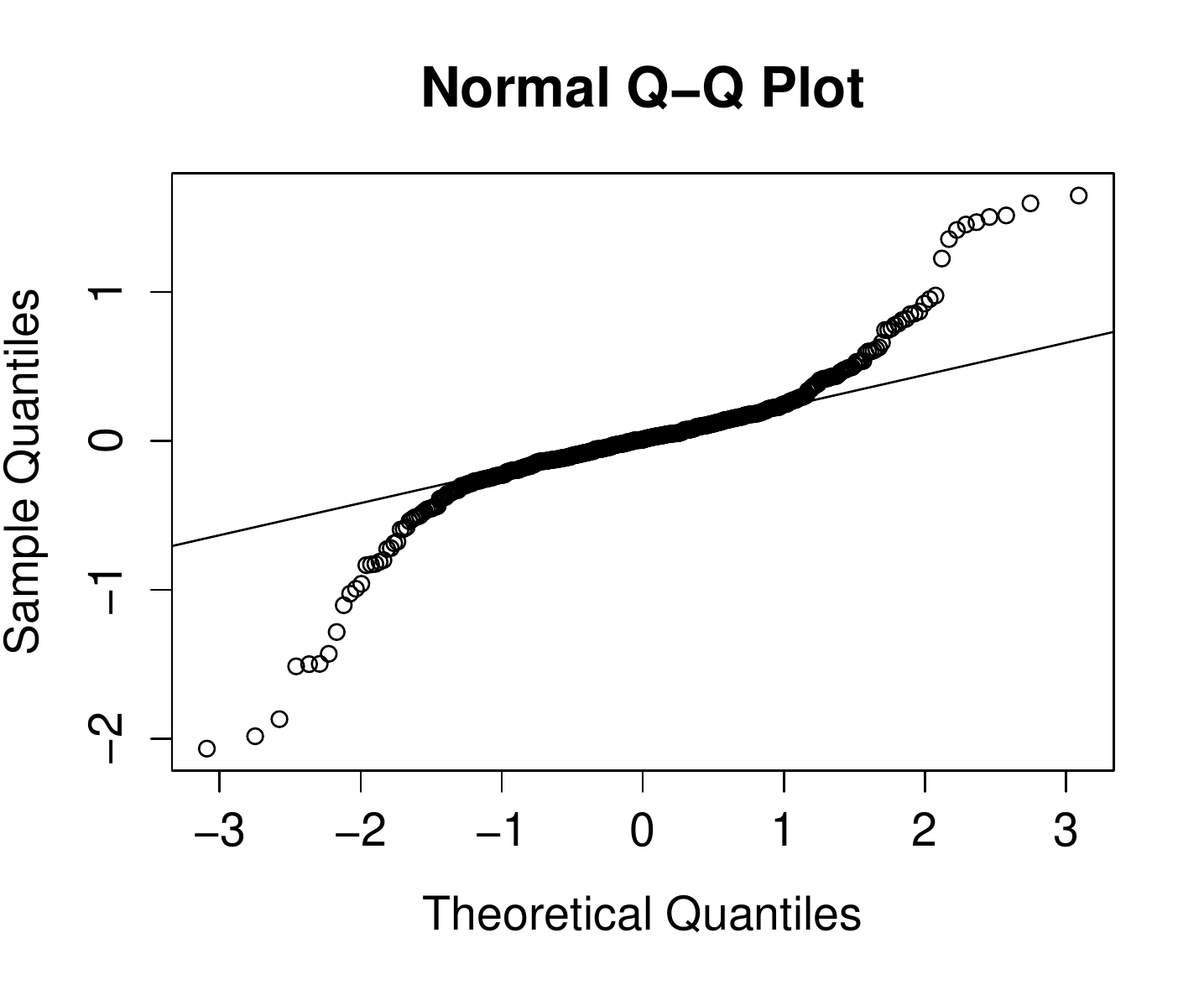}  \end{minipage} \\
   \caption{Q-Q plots of $X_{100,0}(1)$ for Cauchy and Bessel models (from left to right).}\label{fig2}
   \end{figure}

As a simulation example, two types of $\xi(x),$ $x\in \mathbb{R}^2,$ were considered:  long-range dependent Cauchy model with the spectral singularity at zero and cyclical long-range dependent Bessel model with non-zero singularity at 1, see \cite{chi,gne,sch}. 500 realizations of each field were simulated and $X_{100,0}(1)$ were computed for each realisation with the weight function given by (\ref{f0}). We compared empirical distributions of $X_{100,0}(1)$ to the normal law. The observation window was chosen to be large enough ($r=100$) to obtain results close to the asymptotic ones. Figure~\ref{fig2} demonstrates normal Q-Q plots for 500 realisations of each model. Figure~\ref{fig2} suggests that the limit law is normal only when an appropriate weight function is chosen.  The Q-Q plots clearly manifest differences in two types of limit behaviour and support our findings.  
\noindent

\section*{Acknowledgments}
This work was partly supported by La Trobe University
Research Grant "Stochastic Approximation in Finance and Signal Processing".


\end{document}